\documentclass[reqno,12pt]{amsart}

\usepackage{amsmath}
\usepackage{amsfonts}
\usepackage{amssymb}
\usepackage{eufrak}
\usepackage{amscd}
\usepackage{amsthm}
\usepackage{epsfig}
\usepackage{amstext}
\usepackage[all]{xy}

\theoremstyle{plain}

\newtheorem{theorem}{Theorem}
\newtheorem{prop}[theorem]{Proposition}
\newtheorem{cor}[theorem]{Corollary}
\newtheorem{lemma} [theorem] {Lemma}

\newcommand{\qC}{\ensuremath{\rm q\mathbb{C}}}
\newcommand{\qAK}{\ensuremath{ qA\otimes \mathcal K}}
\newcommand{\sdAK}{\ensuremath{\rm S^2A \otimes  \mathcal K}}

\newcommand{\pr}{\left(\begin{smallmatrix} 1& 0
\\ 0 & 0 \end{smallmatrix}\right)}

\newcommand{\C}{\ensuremath{\mathbb C}}

\begin{document}

\title{The $C^*$-algebras $\qAK$ and $\sdAK$ are asymptotically
equivalent}

\author{Tatiana Shulman}

\address{Mathematics Department, New Hampshire University, Durham,
New Hampshire 03824, USA}

\email{tatiana$_-$shulman@yahoo.com}

\subjclass[2000]{46 L80; 19K35}

\keywords {$C^*$-algebra, asymptotic morphism, E-theory, KK-theory}

\date{}

\maketitle

\begin{abstract} Let $A$ be a separable $C^*$-algebra. We prove
that its stabilized second suspension $S^2A\otimes \mathcal K$ and
the $C^*$-algebra $qA\otimes \mathcal K$ constructed by Cuntz in
the framework of his picture of KK-theory are asymptotically
equivalent. This means that there exists an asymptotic morphism
from $S^2A\otimes \mathcal K$ to $qA\otimes \mathcal K$ and an
asymptotic morphism from $qA\otimes \mathcal K$  to $S^2A\otimes
\mathcal K$ whose compositions are homotopic to the identity maps.
This result yields an easy description of the natural
transformation from KK-theory to E-theory. Also by Loring's result
any asymptotic morphism from $\qC$ to any $C^*$-algebra $B$ is
homotopic to a $\ast$-homomorphism. We prove that the same is true
when  $\C$ is replaced by any nuclear $C^*$-algebra $A$ and when
$B$ is stable.
\end{abstract}

\section*{Introduction}
 Let $A$ be a separable $C^*$-algebra. Its first suspension is the
$C^*$-algebra $SA = C_0(\mathbb R)\otimes A$. There are two other
$C^*$-algebras associated to $A$ that are of importance in
KK-theory of Kasparov: the second suspension $C^*$-algebra
$S^2A=C_0(\mathbb R^2)\otimes A$ and the $C^*$-algebra $qA$
constructed by Cuntz \cite{Cuntz} in the framework of his picture
of $KK$-theory. Both $C^*$-algebras can replace $A$ in the
definition of the $KK$-groups: for the second suspension this is
Bott periodicity and for $qA$ this is
 Cuntz's picture for $KK$-theory. These $C^*$-algebras are
$E$-equivalent, i.e. their stabilized suspensions $S^3A\otimes
\mathcal K$ and $SqA\otimes \mathcal K$ are equivalent in the
category of separable $C^*$-algebras with morphisms being homotopy
classes of asymptotic morphisms, where $\mathcal K$ denotes the
$C^*$-algebra of compact operators. In the present paper we show
that they are equivalent in this category without taking the
suspension of the stabilizations. More precisely we construct an
asymptotic morphism from $S^2A\otimes \mathcal K$ to $qA\otimes
\mathcal K$ and a $\ast$-homomorphism from $qA\otimes \mathcal K$
to $S^2A\otimes \mathcal K$ such that their compositions are
homotopic to the identity maps. In general one says that two
$C^*$-algebras are {\it asymptotically equivalent} if there exist
asymptotic morphisms from each to the other whose compositions are
homotopic to the identity maps. So the main result of this paper
(Theorem \ref{main})
 says that $C^*$-algebras $\qAK$ and $\sdAK$
are asymptotically equivalent.

As a corollary  (Corollary \ref{E-theory}) we obtain a description
of E-theory that is similar in form to Cuntz's description of
KK-theory. Cuntz (\cite{Cuntz}) proved that $KK(A, B) = [qA,
B\otimes \mathcal K]$ (where $[\;]$ means  homotopy classes of
$\ast$-homomorphisms). We assert that $E(A, B) = [[qA, B\otimes
\mathcal K]]$ (where $[[\;]]$ means  homotopy classes of asymptotic
morphisms) and that the well known natural transformation $KK(A,
B)\to E(A, B)$ is then nothing but the map that sends any
$\ast$-homomorphism  $qA \to B\otimes \mathcal K$ to itself.

One more corollary (Corollary \ref{Loring}) concerns the question
of when asymptotic morphisms are homotopic to
$\ast$-homomorphisms. In \cite{Lor1} it was proved that any
asymptotic morphism from $\qC$ to any $C^*$-algebra $B$ is
homotopic to a $\ast$-homomorphism. We prove that the same is true
not only for $\C$ but for any nuclear (even K-nuclear)
$C^*$-algebra $A$ if $B$ is assumed to be stable. Recall that a
$C^*$-algebra $B$ is called {\it stable} if $B\otimes \mathcal K
\cong B$.

The plan of the paper is as follows. The first section contains
all necessary information about $C^*$-algebra $qA$. In the second
one we construct an asymptotic morphism $f^A: \sdAK \to \qAK$  and
a $\ast$-homomorphism $g^A: \qAK \to \sdAK$ and show that $f^A$
induces a natural transformation from the $KK$-functor to the
$E$-functor. In the third section we prove that $f^A$ and $g^A$
provide an asymptotic equivalence of the $C^*$-algebras $\sdAK$
and $\qAK$  and obtain the corollaries described above.

\section{Necessary information about $qA$}
 Let $A$ and
$B$ be two $C^*$-algebras. A $C^*$-algebra $C$ is called {\it the
free product of $A$ and $B$} if there are $*$-homomorphisms
$i^A:A\to C$ and
 $i^B:B\to C$ with the following (universal) property: given
$*$-homomorphisms $\phi_A:A\to D$ and $\phi_B:B\to D$ mapping $A$
and $B$ into the same $C^*$-algebra $D$, there is a unique
$*$-homomorphism $\phi: C\to D$ such that $\phi\circ i^A = \phi_A$
and $\phi\circ i^B = \phi_B$.
 The $*$-homomorphisms $i^A$ and $i^B$ are referred to as {\it the
 canonical inclusions}.  The free product of $A$ and $B$ will be denoted
 by $A\ast B$.

Consider $A\ast A$. Let $i_1^A:A\to A\ast A$ and $i_2^A:A\to A\ast
A$ denote the two canonical inclusions of $A$ as a
$C^*$-subalgebra of $A\ast A$. The $C^*$-algebra $qA$ constructed
by Cuntz (\cite{Cuntz}) is the closed ideal in $A\ast A$ generated
by the set $\{i_1(x)-i_2(x): x\in A\}$. One can prove that
elements of the form $$\left(i_1^A(x_1)- i_2^A(x_1)\right) \ldots
\left(i_1^A(x_N)- i_2^A(x_N)\right)$$ and
$$i_1^A(x)\left(i_1^A(x_1)- i_2^A(x_1)\right) \ldots
\left(i_1^A(x_N)- i_2^A(x_N)\right),$$  where $x_0, x_1, \ldots,
x_N\in A$, $N\in \mathbb N$, span a dense $\ast$-subalgebra in
$qA$.

Let $\phi, \psi: A\to B$ be two $*$-homomorphisms. By the
universal property of $A\ast A$ there is a unique $*$-homomorphism
$Q(\phi, \psi): A\ast A\to B$ such that $$Q(\phi, \psi)\circ i_1^A
= \phi,\;\; Q(\phi, \psi)\circ i_2^A = \psi.$$ Let $q(\phi, \psi)$
denote the restriction of $Q(\phi, \psi)$ to $qA$. Note that if
$J$ is an ideal in $B$, then $Q(\phi, \psi)$ maps $qA$ into $J$ if
and only if
 $\phi(x) - \psi(x)\in J$ for all $x\in A$. So in this case, $q(\phi, \psi)\in
Hom(qA, J)$.

\section{Constructing the asymptotic equivalence between $\sdAK$ and $\qAK$}

Below all $C^*$-algebras are assumed to be separable.

For any two $C^*$-algebras $A$ and $B$ Connes and Higson define
$E(A, B)$ to be  the abelian group $[[SA\otimes \mathcal K,
SB\otimes \mathcal K]]$ of homotopy classes of asymptotic morphisms
from $SA\otimes \mathcal K$ to $SB\otimes \mathcal K$ (\cite{CH}).
Recall that an {\it asymptotic morphism from $A$ to $B$} is a family
of maps $(\phi_t)_{t\in [0, \infty)}: A\to B$ satisfying the
following conditions:

\medskip
\hangindent=0.3cm \hangafter=-2 \noindent
 i) for any $a\in A$ the function  $t\mapsto \phi_t(a)$ is continuous;

\medskip
\hangindent=0.3cm \hangafter=-2 \noindent
 ii) for any $a, b\in
A$, $\lambda\in \mathbb C$

\begin{itemize}
\item $\lim_{t\to \infty} \|\phi_t(a^*)-\phi_t(a)^*\| = 0$

\item $\lim_{t\to \infty} \|\phi_t(a + \lambda b)-\phi_t(a) -
\lambda\phi_t(b)\| = 0$;

\item $\lim_{t\to \infty} \|\phi_t(ab)-\phi_t(a)\phi_t(b)\| = 0$.
\end{itemize}

 In \cite{CH} it was also
shown that $[[SA\otimes \mathcal K, SB\otimes \mathcal K]]\cong
[[\sdAK, B\otimes \mathcal K]]$ and we shall always mean  by the
E-group the group $[[\sdAK, B\otimes \mathcal K]]$ of homotopy
classes of asymptotic morphisms from $\sdAK$ to $B\otimes \mathcal
K$.

Let $\beta^{\mathbb C}: C_0(\mathbb R^2)\otimes \mathcal K\to
\mathcal K$ be the Bott asymptotic morphism. In fact it is the
tensor product of the identity map $id_{\mathcal K}: \mathcal K\to
\mathcal K$ with the restriction to $C_0(\mathbb R^2) \subset
C(\mathbb T^2)$ of the family of maps from $C(\mathbb T^2)$ to
$\mathcal K^+$ constructed in the Voiculescu's example of almost
commuting unitaries (\cite{Voic}), but here we shall not use an
explicit form of $\beta^{\C}$ but only the fact that it induces
the identity map in the K-groups. Let
 $$\beta ^A = \beta^{\mathbb C}\otimes
id_A: S^2A\otimes \mathcal K \to A\otimes \mathcal K.$$ Obviously
$\beta^A\in E(A, A)$. Note that since we always consider asymptotic
morphisms up to homotopy we denote in the same way a class of
homotopy equivalent asymptotic morphisms and any its representative.


For the KK-groups we will use Cuntz's approach (\cite{Cuntz}) in
which, as already was written, one regards   $KK(A, B)$ as the
group $[qA\otimes \mathcal K, B\otimes \mathcal K]$ of homotopy
classes of $*$-homomorphisms from $qA\otimes \mathcal K$ to
$B\otimes \mathcal K$. Let
$$\gamma^A = q(id_A, 0)\otimes id_{\mathcal K}: \qAK \to A\otimes \mathcal K.$$
Then $\gamma^A\in KK(A, A)$ and it is a unit element for the
associative product $KK(A, B)\times KK(B, C)\to KK(A, C)$. Namely
there exists a bilinear pairing $KK(A, B)\times KK(B, C)\to KK(A,
C)$ such that $x\times \gamma^B = x = \gamma^A\times x$ for any
$x\in KK(A, B)$ (\cite{Cuntz}).

 Let $A$ be a $C^*$-algebra.  By \cite{CH}
there exists a natural transformation from the functor $KK(A,-)$
into the functor $E(A, -)$ which is unique up to its value on
 $\gamma^A\in KK(A, A)$.
Let $$I_{A, B}:KK(A, B)\to E(A, B)$$ be such a natural
transformation that $I_{A, A}(\gamma^A) = \beta^A.$ Define an
asymptotic morphism $f^A:\sdAK\to \qAK$ by $$f^A = I_{A,
qA}(id_{\qAK}).$$

\noindent The following easy theorem asserts that the  asymptotic
morphism $f^A$ induces the natural transformation $I_{A, B}$.

\begin{theorem}\label{ind} $I_{A, B}(\phi) = \phi\circ f^A$ for any $\phi \in
KK(A, B)$.
\end{theorem}

{\bf Proof}. Since $\phi\in KK(A, B)$ is a $\ast$-homomorphism
from $\qAK$ to $B\otimes \mathcal K$  it induces the maps
$\phi_{KK}:KK(A, qA)\to KK(A, B)$ and $\phi_E: E(A, qA)\to E(A,
B)$ in the KK-groups and the E-groups respectively. By the
definition of a natural transformation of covariant functors the
following diagram commutes
$$\begin{CD} KK(A, B) @>{I_{A, B}}>> E(A, B) \\
@A{\phi_{KK}}AA  @AA{\phi_E}A \\ KK(A, qA) @>{I_{A, qA}}>> E(A,
qA)
\end{CD} $$

\noindent Hence for the element $id_{\qAK}\in KK(A, qA)$ we get
$$\phi_E(I_{A, qA}(id_{\qAK})) = I_{A, B}(\phi_{KK}(id_{\qAK}).$$

\noindent But $\phi_E(I_{A, qA}(id_{\qAK})) = \phi \circ I_{A,
qA}(id_{\qAK}) = \phi\circ f^A$ and $I_{A, B}(\phi_{KK}(id_{\qAK})
= I_{A, B}(\phi\circ id_{\qAK}) = I_{A, B}(\phi)$. \qed

\begin{cor}\label{sure} $\gamma^A\circ f^A = \beta^A$.
\end{cor}
{\bf Proof.} By Theorem \ref{ind} $\gamma^A\circ f^A = I_{A,
A}(\gamma^A)$. Since we have chosen a natural transformation to be
equal $\beta^A$ on the element $\gamma^A$ we get $\gamma^A\circ
f^A = \beta^A$. \qed

\bigskip
Now we define a $\ast$-homomorphism $g^A:\qAK\to \sdAK$ in the
following way.
 Let $\pi_1, \pi_2: \C \to C_0(\mathbb R^2)^+\otimes M_2$ be two
$*$-homomorphisms given by
$$\pi_1(1) = \left(\begin{array}{cc} 1& 0
\\ 0 & 0 \end{array}\right),\;\;\; \pi_2(1) = p_{Bott} =
\frac{1}{1+z\bar{z}}\left(\begin{array}{cc} z\bar{z} & z \\
\bar{z} & 1\end{array}\right)$$ (we identify $\mathbb R^2$ with
$\mathbb C$). Fix once and for all some inclusion $j: M_2 \to
\mathcal K$ and some isomorphism $i:\mathcal K\otimes \mathcal K\to
\mathcal K$. Define $\tilde\pi_1, \tilde\pi_2:A \to A\otimes
C_0(\mathbb R^2)^+\otimes \mathcal K$ by
$$\tilde \pi_1 = (j\otimes id_{A\otimes C_0(\mathbb
R^2)^+})\circ(id_A\otimes \pi_1), \;\;\; \tilde \pi_2 = (j\otimes
id_{A\otimes C_0(\mathbb R^2)^+})\circ(id_A\otimes \pi_2)$$
respectively.
   Since $$\tilde\pi_1(a) - \tilde\pi_2(a) \in C_0(\mathbb R^2)\otimes \mathcal K\otimes A = \sdAK$$ for any $a\in A$, the $\ast$-homomorphism
$q(\tilde\pi_1, \tilde\pi_2): qA \to S^2A\otimes \mathcal K$ is
defined.

\noindent Set
 $$g^A =(id_{S^2A}\otimes i) \circ
(q(\tilde\pi_1, \tilde\pi_2)\otimes id_{\mathcal K}).$$

\medskip

\noindent In the next section we show that $ f^A $ and $g^A$ provide
an asymptotic equivalence between $\sdAK$ and $qA\otimes \mathcal
K$.

\section{Proof of the main assertion}
To prove that $ f^A $ and $g^A$ provide an asymptotic equivalence
between $\sdAK$ and $qA\otimes \mathcal K$ we are going to show
that their compositions induce the identity maps in $E$-functor
and in the functor $G$ that will be introduced in subsection
\ref{subsG}.
\subsection{The maps induced by $f^A$ and $g^A$ in E-functor}
\begin{lemma}\label{K-theory} $\beta^A\circ g^A \sim \gamma^A$.
\end{lemma}
{\bf Proof.} Note first of all that $g^A:\qAK \to \sdAK$ and
$\gamma^A:\qAK\to A\otimes \mathcal K$ factorize through the
$C^*$-algebra $q\mathbb C\otimes A\otimes \mathcal K$. Namely let
$\eta_1, \eta_2: A\to (\mathbb C\ast \mathbb C)\otimes A$ be given
by formulas $$\eta_1(a) = i_1^{\mathbb C}(1)\otimes a,\;\; \eta_2(a)
= i_2^{\mathbb C}(1)\otimes a$$ for any $a\in A$.  Set
$$s^A = q(\eta_1, \eta_2): qA\to q\mathbb C\otimes A.$$
It is easy to see that the diagrams $$\xymatrix{qA\otimes
K\ar[rr]^-{\gamma^A}\ar[dr]_-{s^A\otimes id_K} & &A\otimes K\\
&q\C\otimes A\otimes K\ar[ru]_-{\gamma^{\C}\otimes id_A}}$$ and
$$\xymatrix{qA\otimes K\ar[rr]^-{g^A}\ar[dr]_-{s^A\otimes id_K} &
& S^2A\otimes K\\&q\C\otimes A\otimes K\ar[ur]_-{g^{\C}\otimes
id_A}}$$ commute, that is  $$\gamma^A = (\gamma^{\mathbb C}\otimes
id_A)\circ (s^A\otimes id_{\mathcal K}), \;\;\;\;\;g^A =
(g^{\mathbb C}\otimes id_A)\circ (s^A\otimes id_{\mathcal K}).$$
Since $\beta^A = \beta^{\mathbb C}\otimes id_A$ we have to
establish the homotopy equivalence
$$(\gamma^{\mathbb
C}\otimes id_A)\circ (s^A\otimes id_{\mathcal K}) \sim
(\beta^{\mathbb C}\otimes id_A)\circ (g^{\mathbb C}\otimes
id_A)\circ (s^A\otimes id_{\mathcal K})$$ or, equivalently,
$$\gamma^{\mathbb C}\sim \beta^{\C}\circ g^{\C}.$$ For that we use
K-theory. Let $\gamma^{\mathbb C}_{\ast}$ and $(\beta^{\C}\circ
g^{\C})_{\ast}$ be the induced homomorphisms from $K_0(q\C)$ to
$K_0(\C)$. For the generator $[i_1^{\C}(1)] - [i_2^{\C}(1)]$ of
$K_0(q\C)$ we have $$ (\beta^{\C}\circ g^{\C})_{\ast}([i_1^{\C}(1)]
- [i_2^{\C}(1)]) = \beta^{\C}_{\ast}([\pr] -[p_{Bott}]) = [1],$$
$$\gamma^{\mathbb C}_{\ast}([i_1^{\C}(1)] - [i_2^{\C}(1)]) = [1]
-[0] = [1].$$ We used here that $[\pr] -[p_{Bott}]$ is a generator
of $K_0(S^2\mathbb C)$ and that the Bott map $\beta^{\C}$ induces
the identity homomorphism in K-theory. So $\gamma^{\mathbb C}$ and
$\beta^{\C}\circ g^{\C}$ induce the same homomorphisms in
K-theory. This implies that these asymptotic homomorphisms are
homotopic because, by Universal coefficients theorem,
$$Hom(K_0(q\C), K_0(\mathcal K))\oplus Hom(K_1(q\C), K_1(\mathcal K)) \cong KK(q\C,
\mathcal K)\oplus KK(Sq\C, \mathcal K),$$ and since $$K_1(q\C) =
K_1(\mathcal K) = 0,\; KK(Sq\C, \mathcal K) = 0,\; KK(q\C, \mathcal
K) = [q\C, \mathcal K] \stackrel{[3]}{=}
 [[q\C, \mathcal
K]]$$ we get
$$Hom(K_0(q\C), K_0(\C)) \cong [[q\C\otimes \mathcal K, \mathcal K]].$$ \qed

\medskip

Let $B$ be any $C^*$-algebra. Let $f^A_E: E(B, S^2A) \to E(B, qA)$
and $g^A_E: E(B, qA) \to E(B, S^2A)$ be the maps induced by $f^A$
and $g^A$ respectively.

\begin{prop}\label{E}  $f^A_E\circ g^A_E = id$,

$\;\;\;\;\;\;\;\;\;\;\;\;\;\;\;\;\;\;\;\;\; g^A_E\circ f^A_E =
id$.
\end{prop}

Here $id$ means both the identity map from $E(B, S^2A)$ into
itself and the identity map from $E(B, qA)$ into itself.

{\bf Proof.} Consider the following diagram

$$
\xymatrix{ & E(B, A) & \\
 E(B, S^2A)\ar[ur]^{\beta^A_E} \ar@<-1ex>[rr]_{f_E^A}& &E(B,
qA) \ar[ul]_{\gamma_E^A} \ar@<-1ex>[ll]_{g_E^A}\\
 }
$$

\noindent Here $\beta_E^A$ and $\gamma_E^A$ are the maps induced by
$\beta^A$ and $\gamma^A$ respectively. It is proved in \cite{CH}
that $\beta_E^A$ is an isomorphism. Furthermore  $\gamma_E^A$ also
is an isomorphism. Indeed by \cite{Cuntz} the map induced by
$\gamma^A$ in any covariant, homotopy invariant, split exact and
stable functor is an isomorphism. Since the functor $E(B, -)$ has
all these properties $\gamma_E^A$ is an isomorphism.

\noindent  By Lemma \ref{K-theory}, $\beta^A_E\circ g^A_E =
\gamma^A_E$ whence \begin{equation}\label{old1}g^A_E =
(\beta_E^A)^{-1}\circ \gamma_E^A \end{equation} By Corollary
\ref{sure} $\gamma^A_E\circ f^A_E = \beta^A_E$ whence
\begin{equation}\label{old2}f^A_E = (\gamma^A_E)^{-1}\circ \beta_E^A
\end{equation}

\noindent The assertions of the proposition follow from (\ref{old1})
and (\ref{old2}). \qed
\subsection{The maps induced by $f^A$ and $g^A$ in G-functor}\label{subsG}
Now instead of E-functor we are going to consider another
bifunctor $G(B, A)$   and prove the result similar to Lemma
\ref{E} for the maps, induced by $f^A$ and $g^A$ in the functor
$G(B, -)$, where $B$ is  fixed. Namely let $G(B, A)$ be the
semigroup $[[qB\otimes \mathcal K, A\otimes \mathcal K]]$ of the
classes of homotopy equivalent asymptotic homomorphisms from
$qB\otimes \mathcal K$ to $A\otimes \mathcal K$. Obviously this is
a contravariant functor in the first variable and a covariant
functor in the second one. We need two results about this
bifunctor
--- the Bott periodicity and the isomorphism $G(B, A)\cong G(B, qA)$. To
prove them we need first of all a construction which produces an
asymptotic morphism $q\psi: qD_1 \to qD_2$ out of an asymptotic
morphism $\psi:D_1 \to D_2$, where $D_1, D_2$ are any
$C^*$-algebras.

An asymptotic morphism $\psi$ gives rise to a genuine
$\ast$-homomorphism $$F: D_1 \to C_b([0, \infty), D_2)/C_0([0,
\infty), D_2)$$ given by $$F(x) = \psi_t(x) + C_0([0, \infty),
D_2)$$ for any $x\in D_1$. There are two $\ast$-homomorphisms
$\overline{i_1}, \overline{i_2}: C_b([0, \infty), D_2)\to C_b([0,
\infty), D_2\ast D_2)$ given by formulas
$$\overline{i_1}(f)(t) = i_1^{D_2}(f(t)),\;\; \overline{i_2}(f)(t)
= i_2^{D_2}(f(t)),$$ $f\in C_b([0, \infty), D_2)$.  Since these
$\ast$-homomorphisms send $C_0([0, \infty), D_2)$ to $C_0([0,
\infty), D_2\ast D_2)$ we have two $\ast$-homomorphisms
$$\hat{i_1}, \hat{i_2}: C_b([0, \infty), D_2)/C_0([0, \infty), D_2)
\to C_b([0, \infty), D_2\ast D_2)/C_0([0, \infty), D_2\ast D_2).$$
Set $$\Phi = Q(\hat{i_1}\circ F, \hat{i_2}\circ F): D_1\ast D_1 \to
C_b([0, \infty), D_2\ast D_2)/C_0([0, \infty), D_2\ast D_2).$$ Let
$p:C_b([0, \infty), D_2\ast D_2) \to C_b([0, \infty), D_2\ast
D_2)/C_0([0, \infty), D_2\ast D_2)$ be the canonical surjection.
Since $$\Phi(i_1^{D_1}(a)) = p\left(i_1^{D_2}(\psi_t(a))\right),\;\;
 \Phi(i_2^{D_1}(a)) = p\left(i_2^{D_2}(\psi_t(a))\right)$$  for any $a\in D_1$,
 and since $qD_1$ is the closed ideal generated by  the set $\{i_1^{D_1}(a)-i_2^{D_1}(a):
a\in D_1\}$,   we get $$\Phi(qD_1) \subset p\left(C_b([0, \infty),
qD_2)\right).$$ We shall denote the restriction of $\Phi$ to
$qD_1$ also by $\Phi$. Define a $\ast$-homomorphism
$$\tau: p\left(C_b([0, \infty), qD_2)\right)\to C_b([0, \infty),
qD_2)/C_0([0, \infty), qD_2)$$ by
$$\tau(p(f)) = f + C_0([0, \infty), qD_2),$$ $f\in C_b([0, \infty),
qD_2)$. It is well-defined because for any  $f\in C_b([0, \infty),
qD_2)$ the condition $f\in C_0([0, \infty), D_2\ast D_2)$ implies
$f\in C_0([0, \infty), qD_2)$. So we have  $\tau\circ \Phi: qD_1 \to
C_b([0, \infty), qD_2)/C_0([0, \infty), qD_2)$.
  Choose a continuous section
$$s: C_b([0, \infty), qD_2)/C_0([0, \infty), qD_2) \to C_b([0,
\infty), qD_2)$$ (it exists by Bartle-Graves theorem, \cite{Bartle,
Lor3}) and define an asymptotic morphism $q\psi$ by
$$(q\psi)_t(x) = \bigl(s\left(\tau\circ \Phi(x)\right)\bigl)(t).$$
Thus we get an asymptotic morphism $q\psi:qD_1 \to qD_2$ out of an
asymptotic morphism $\psi:D_1 \to D_2$.

For any $C^*$-algebra $D\;$ let $$\rho^D = q( i_1^D\otimes
id_{\mathcal K}, i_2^D\otimes id_{\mathcal K} ): q(D\otimes
\mathcal K) \to qD\otimes \mathcal K$$ and let $\theta_D:
qD\otimes \mathcal K \to q^2D\otimes \mathcal K$ denote the
isomorphism constructed in \cite{Cuntz}.

\begin{lemma}\label{1} The diagram $$\xymatrix{q(A\otimes
K)\otimes K\ar[r]^-{\gamma^{A\otimes K}}\ar[d]^-{\rho^A\otimes
id_K}&A\otimes K\otimes K\ar[r]^-{id_A\otimes i}&A\otimes K\\
qA\otimes K \otimes K\ar[d]^{id_{qA}\otimes i}&&\\
qA\otimes K\ar[uurr]_{\gamma^A}&& }$$
 is commutative, namely $\gamma^A\circ (id_{qA}\otimes i)\circ (\rho^A\otimes
id_{\mathcal K}) = (id_A\otimes i)\circ \gamma^{A\otimes \mathcal
K}$.
\end{lemma}
{\bf Proof.} Since elements of the form $$\left((i_1^{A\otimes
\mathcal K}(a\otimes T) - i_2^{A\otimes \mathcal K}(a\otimes
T)\right)\otimes S$$ and $$\left(i_1^{A\otimes \mathcal
K}(a_0\!\otimes \!T_0)\left(i_1^{A\otimes \mathcal K}(a\!\otimes
\!T) - i_2^{A\otimes \mathcal K}(a\!\otimes
\!T)\right)\right)\otimes S ,$$ where $T, S, T_0\in \mathcal K$,
$a, a_0\in A$, span a dense subspace of $q(A\otimes \mathcal
K)\otimes \mathcal K$ (see \cite{Th}, for example) it is enough to
check that $\gamma^A\circ (id_{qA}\otimes i)\circ (\rho^A\otimes
id_{\mathcal K})$ and $(id_A\otimes i)\circ \gamma^{A\otimes
\mathcal K}$ coincide on elements of such form. For any $T, S\in
\mathcal K$, $a\in A$ we have
\begin{multline*} \gamma^A\circ (id_{qA}\otimes i)\circ
(\rho^A\otimes id_{\mathcal K})\Bigl(\left(i_1^{A\otimes \mathcal
K}(a\otimes
T) - i_2^{A\otimes \mathcal K}(a\otimes T)\right)\otimes S \Bigl) =\\
a\otimes i(T\otimes S) =\\ (id_A\otimes i)\circ \gamma^{A\otimes
\mathcal K}\Bigl(\left(i_1^{A\otimes \mathcal K}(a\otimes T) -
i_2^{A\otimes \mathcal K}(a\otimes T)\right)\otimes S \Bigl),
\end{multline*} for another pair $T_0\in \mathcal K$, $a_0\in A$ we have  \begin{multline*} \gamma^A\circ (id_{qA}\otimes i)\circ
(\rho^A\otimes id_{\mathcal K})\Bigl(\!\left(i_1^{A\otimes
\mathcal K}(a_0\!\otimes \!T_0)\left(i_1^{A\otimes \mathcal
K}(a\!\otimes \!T) - i_2^{A\otimes \mathcal K}(a\!\otimes
\!T)\right)\right)\otimes S \!\Bigl) = \\ a_0a\otimes
i(T_0T\otimes S) =\\ (id_A\otimes i)\circ \gamma^{A\otimes
\mathcal K}\Bigl(\left(i_1^{A\otimes \mathcal K}(a_0\otimes
T_0)\left(i_1^{A\otimes \mathcal K}(a\otimes T) - i_2^{A\otimes
\mathcal K}(a\otimes T)\right)\right)\otimes S \Bigl)
\end{multline*} and we are done.\qed

\begin{lemma}\label{2} Let $\phi \in [[qB, A\otimes \mathcal K]]$.
Then the diagram $$\xymatrix {q^2B\otimes K \ar[rr]^-{q\phi\otimes
id_K} \ar[d]^{\gamma^{qB}} && q(A\otimes K)\otimes K
\ar[rr]^{\gamma^{A\otimes K}} && A\otimes K\otimes K \\
qB\otimes K \ar@/_/[urrrr]_{\phi\otimes id_K}&& & &}$$
 commutes, that is $\gamma^{A\otimes \mathcal
K}\circ (q\phi \otimes id_{\mathcal K}) = (\phi\otimes
id_{\mathcal K})\circ \gamma^{qB}$.
\end{lemma}
{\bf Proof.} Let $x\in qB$, $T\in \mathcal K$. By the definition
of $q\phi$ we have $$(q\phi)_t\left(i_1^{qB}(x)-i_2^{qB}(x)\right)
- \left(i_1^{A\otimes \mathcal K}(\phi_t(x))-i_2^{A\otimes
\mathcal K}(\phi_t(x))\right)\to 0$$ when $t\to \infty$. Hence
\begin{multline*} \lim_{t\to \infty} \Bigl[\gamma^{A\otimes \mathcal K}\circ ((q\phi)_t \otimes
id_{\mathcal K})\left(\left(i_1^{qB}(x)-i_2^{qB}(x)\right)\otimes T\right) -\\
(\phi_t\otimes id_{\mathcal K})\circ
\gamma^{qB}\left(\left(i_1^{qB}(x)-i_2^{qB}(x)\right)\otimes
T\right)\Bigl] =\\ \lim_{t\to \infty}\Bigl[ \gamma^{A\otimes
\mathcal K}\Bigl(\left(i_1^{A\otimes \mathcal
K}(\phi_t(x))-i_2^{A\otimes \mathcal K}(\phi_t(x))\right)\otimes
T\Bigl) - \phi_t(x)\otimes T\Bigl] = 0.
\end{multline*} In a similar way we find that $\gamma^{A\otimes \mathcal K}\circ (q\phi \otimes id_{\mathcal K})$ and $(\phi\otimes id_{\mathcal K})\circ
\gamma^{qB}$ asymptotically agree on elements
$\left(i_1^{qB}(x_0)\left(i_1^{qB}(x)-i_2^{qB}(x)\right)\right)\otimes
T$ when $x_0, x\in qB$, $T\in \mathcal K$. Since elements of the
form $(\left(i_1^{qB}(x)-i_2^{qB}(x)\right)\otimes T$ and
$\left(i_1^{qB}(x_0)\left(i_1^{qB}(x)-i_2^{qB}(x)\right)\right)\otimes
T$ span a dence subspace of $qB\otimes \mathcal K$ we see that
$\gamma^{A\otimes \mathcal K}\circ (q\phi \otimes id_{\mathcal K})
= (\phi\otimes id_{\mathcal K})\circ \gamma^{qB}$. \qed


\begin{lemma}\label{4} Let $\phi \in [[qB, qA\otimes \mathcal K]]$.
Then the diagram
$$\xymatrix{q^2B\ar[r]^-{q\phi}\ar[d]^{q(id_{qB},\; 0)} &
q(qA\otimes K)\ar[r]^{\rho^{qA}}&q^2A\otimes
K\ar[r]^{\gamma^{qA}}&qA\otimes K\\
qB\ar@/_/[urrr]_-{\phi}&&&}$$ is commutative, that is
 $\gamma^{qA}\circ\rho^{qA}\circ
q\phi = \phi\circ q(id_{qB}, 0)$.
\end{lemma}
{\bf Proof.} Let $x\in qB$, $t\in [0, \infty)$. Writing
$\phi_t(x)$ in the form $$\phi_t(x) = \lim_{k\to \infty}
\sum_{i=1}^{N_k} z_i^{(k)}(t)\otimes T_i^{(k)}(t),$$ where
$z_i^{(k)}(t)\in qA$, $T_i^{(k)}(t)\in \mathcal K$, we get
\begin{multline*} \lim_{t\to \infty}\Bigl[\gamma^{qA}\circ\rho^{qA}\circ (q\phi)_t\left(i_1^{qB}(x)-
i_2^{qB}(x)\right) - \phi_t\circ q(id_{qB}, 0)\left(i_1^{qB}(x)-
i_2^{qB}(x)\right)\Bigl] =\\ \lim_{t\to \infty}\Bigl[
\gamma^{qA}\circ\rho^{qA}\left(i_1^{\qAK}(\phi_t(x)) -
i_2^{\qAK}(\phi_t(x))\right) - \phi_t(x)\Bigl] = \\\lim_{t\to
\infty}\Bigl[ \lim_{k\to \infty}\sum_{i=1}^{N_k}\gamma^{qA}
\left(i_1^{qA}(z_i^{(k)})\!\otimes \!T_i^{(k)} -
i_2^{qA}(z_i^{(k)})\!\otimes \!T_i^{(k)}\right) - \lim_{k\to
\infty} \sum_{i=1}^{N_k} z_i^{(k)}\otimes T_i^{(k)}\Bigl] =
0.\end{multline*} In a similar way we find that
$\gamma^{qA}\circ\rho^{qA}\circ q\phi$ and $\phi\circ q(id_{qB},
0)$ asymptotically agree on elements
$i_1^{qB}(x_0)\left(i_1^{qB}(x)- i_2^{qB}(x)\right)$, where $x_0,
x\in qB$. Since elements of the form $i_1^{qB}(x)- i_2^{qB}(x)$
and $i_1^{qB}(x_0)\left(i_1^{qB}(x)- i_2^{qB}(x)\right)$ span a
dense subspace of $qB$ we conclude that the asymptotic morphisms
$\gamma^{qA}\circ\rho^{qA}\circ q\phi$ and $\phi\circ q(id_{qB},
0)$ coincide. \qed

\bigskip

Let $\psi\in G(B, A)$. There is an asymptotic morphism $\phi:qB\to
A\otimes \mathcal K$ such that $(id_A\otimes i)\circ(\phi\otimes
id_{\mathcal K})\sim \psi$ (\cite{CH}). Define an asymptotic
morphism $\Gamma(\psi)\in G(B, qA)$ by the following composition
$$\xymatrix @1{qB\otimes K\ar[r]^-{\theta_B} & q^2B\otimes
K\ar[rr]^-{q\phi\otimes id_K}&&q(A\otimes K)\otimes
K\ar[rr]^-{\rho^A\otimes id_K} && qA\otimes K\otimes
K\ar[rr]^-{id_{qA}\otimes i}& &qA\otimes K}.$$ Thus a map
$\Gamma: G(B, A)\to G(B, qA)$ is defined by formula
$$ \Gamma(\psi) = (id_{qA}\otimes i)\circ (\rho^A\otimes
id_{\mathcal K})\circ (q\phi\otimes id_{\mathcal K})\circ
\theta_B$$ for any $\psi\in G(B, A)$.

\noindent Let $\gamma^A_G: G(B, qA)\to G(B, A)$ be the map induced
by $\gamma^A$.

\begin{prop}\label{Gq}  $\Gamma: G(B, A) \to G(B, qA)$ is a semigroup
isomorphism with inverse $\gamma^A_G$.
\end{prop}
{\bf Proof.} Obviously $\Gamma$ and $\gamma^A_G$ are semigroup
homomorphisms so we have to check only the following:

\medskip

  (i) $\Gamma(\gamma^A_G(\psi))\sim \psi$ for any $\psi\in
G(B, qA)$,
\medskip

 (ii) $\gamma^A_G(\Gamma(\psi))\sim \psi$ for any $\psi\in G(B,
A)$.

\medskip

\noindent (i): Let $\psi\in G(B, qA)$ and $\phi:qB\to \qAK$ be such
an asymptotic morphism that $(id_{qA}\otimes i)\circ(\phi\otimes
id_{\mathcal K})\sim \psi$. Then
\begin{multline*}
\Gamma(\gamma^A_G(\psi)) = (id_{qA}\otimes i)\circ (\rho^A\otimes
id_{\mathcal K})\circ (q(\gamma^A\circ \phi)\otimes id_{\mathcal K})\circ \theta_B = \\
(id_{qA}\otimes i)\circ (\rho^A\otimes id_{\mathcal K})\circ
(q\gamma^A\otimes id_{\mathcal K})\circ (q\phi\otimes id_{\mathcal
K})\circ \theta_B,\end{multline*} because clearly $q(\gamma^A\circ
\phi)\otimes id_K =  (q\gamma^A\otimes id_K)\circ (q\phi\otimes
id_K)$.

By (\cite{Th}, Lemma 5.1.11) $\rho^A\circ q\gamma^A \sim
\gamma^{qA}\circ\rho^{qA}$ and we have
\begin{multline*}\Gamma(\gamma^A_G(\psi)) =  (id_{qA}\otimes i)\circ (\gamma^{qA}\otimes id_{\mathcal K})\circ
(\rho^{qA}\otimes id_{\mathcal K})\circ
(q\phi\otimes id_{\mathcal K})\circ \theta_B \stackrel{Lemma \ref{4}} = \\
(id_{qA}\otimes i)\circ (\phi\otimes id_{\mathcal K})\circ
\gamma^{qB}\circ \theta_B \sim  \psi\circ \gamma^{qB}\circ \theta_B
\stackrel{[1]} {\sim} \psi.
\end{multline*}
\noindent (ii): Now let $\psi\in G(B, A)$ and $\phi:qB\to A\otimes
\mathcal K$ be such an asymptotic morphism that $(id_A\otimes
i)\circ (\phi\otimes id_{\mathcal K}) \sim \psi$. Then
\begin{multline*} \gamma^A_G(\Gamma(\psi)) =  \gamma^A\circ
(id_{qA}\otimes i)\circ (\rho^A\otimes id_{\mathcal K})\circ
(q\phi\otimes id_{\mathcal K})\circ \theta_B \stackrel{Lemma \ref{1}} =  \\
(id_A\otimes i)\circ \gamma^{A\otimes \mathcal K} \circ
(q\phi\otimes id_{\mathcal K}) \circ
 \theta_B \stackrel{Lemma \ref{2}} = \\(id_A\otimes i) \circ
(\phi\otimes id_{\mathcal K}) \circ \gamma^{qB} \circ \theta_B \sim
\psi \circ \gamma^{qB }\circ \theta_B \stackrel{[1]}{\sim} \psi.
\end{multline*} \qed

\begin{lemma}\label{new} The diagram $$\xymatrix{q(S^2A\otimes
K)\otimes K\ar[rr]^-{q\beta^A\otimes
id_K}\ar[dd]_-{\gamma^{S^2A\otimes K}} & & q(A\otimes K)\otimes K
\ar[rr]^-{\rho^A\otimes id_K}& & qA\otimes K\otimes K
\ar[d]^-{id_{qA\otimes i}}\\ && &&qA\otimes K\ar[d]^-{S^2A\otimes
K}\\
 S^2A\otimes K\otimes K\ar[rrrr]^-{id_{S^2A}\otimes i} & & & &S^2A\otimes
K }$$ commutes.

\noindent Namely
 $g^A\circ (id_{qA}\otimes i)\circ
(\rho^{A}\otimes id_{\mathcal K})\circ (q\beta^A\otimes id_{\mathcal
K}) \sim (id_{S^2A}\otimes i)\circ \gamma^{S^2A\otimes \mathcal K}$.
\end{lemma}
{\bf Proof.} We will prove the assertion by establishing the
commutativity of the left and right tiangles of the diagram
$$\xymatrix{q(S^2A\otimes K)\otimes K\ar[rr]^-{q\beta^A\otimes
id_K}\ar[dd]_-{\gamma^{S^2A\otimes K}} & & q(A\otimes K)\otimes K
\ar[rr]^-{\rho^A\otimes id_K}& & qA\otimes K\otimes K
\ar[d]^-{id_{qA\otimes i}}\ar@{-->}[ddllll]_-{g^A\otimes id_K}\\
&& &&qA\otimes K\ar[d]^-{S^2A\otimes
K}\\
 S^2A\otimes K\otimes K\ar[rrrr]^-{id_{S^2A}\otimes i} & & & &S^2A\otimes
K. }$$ To prove the commutativity of the right triangle we have to
prove
\begin{equation}\label{n_0} g^A\circ (id_{qA}\otimes i)\sim
(id_{S^2A}\otimes i)\circ (g^A\otimes id_{\mathcal
K}).\end{equation} Let $h_1, h_2: K\otimes \mathcal K\otimes
\mathcal K \to \mathcal K$ be the isomorphisms which send
$T_1\otimes T_2\otimes T_3$ to $i(T_1\otimes i(T_2\otimes T_3))$
and $i(i(A\otimes B)\otimes C)$ respectively for any operators
$T_1, T_2, T_3\in \mathcal K.$
 Then for any
$T, S\in \mathcal K$, $a\in A$ we have
\begin{multline*}\left(id_{S^2A}\otimes \left(h_2\circ h_1^{-1}\right)\right)\circ g^A\circ
(id_{qA}\otimes i)\Bigl(\left(i_1^A(a)-i_2^A(a)\right)\otimes
T\otimes S\Bigl) = \\\left(id_{S^2A}\otimes \left(h_2\circ
h_1^{-1}\right)\right)\Bigl(a\otimes i\Bigl(j\pr\otimes i(T\otimes
S)\Bigl)\Bigl) = \\ a\otimes i\Bigl(i(j\pr \otimes T)\otimes S\Bigl)
= (id_{S^2A}\otimes i)\circ (g^A\otimes id_{\mathcal
K})\Bigl(\left(i_1^A(a)-i_2^A(a)\right)\otimes T\otimes
S\Bigl)\end{multline*} and for another $a_0\in A$
\begin{multline*}\left(id_{S^2A}\otimes \left(h_2\circ
h_1^{-1}\right)\right)\circ g^A\circ (id_{qA}\otimes
i)\Bigl(i_1^A(a_0)\left(i_1^A(a)-i_2^A(a)\right)\otimes T\otimes
S\Bigl) = \\\left(id_{S^2A}\otimes \left(h_2\circ
h_1^{-1}\right)\right)\Bigl(a_0a\otimes i\Bigl(j\bigl(\pr(\pr -
p_{Bott})\bigl)\otimes i(T\otimes S)\Bigl)\Bigl) = \\ a_0a\otimes
i\Bigl(i\left(j\bigl(\pr(\pr - p_{Bott})\bigl) \otimes
T\right)\otimes S\Bigl)= \\ (id_{S^2A}\otimes i)\circ (g^A\otimes
id_{\mathcal
K})\Bigl(i_1^A(a_0)\left(i_1^A(a)-i_2^A(a)\right)\otimes T\otimes
S\Bigl).\end{multline*}
 Since elements of the form
$$\left(i_1^A(a)-i_2^A(a)\right)\otimes T\otimes S$$ and
$$i_1^A(a_0)\left(i_1^A(a)-i_2^A(a)\right)\otimes
T\otimes S$$ span a dense subspace of $qA\otimes \mathcal K\otimes
\mathcal K$ we get
$$\left(id_{S^2A}\otimes \left(h_2\circ
h_1^{-1}\right)\right)\circ g^A\circ (id_{qA}\otimes i) =
(id_{S^2A}\otimes i)\circ (g^A\otimes id_{\mathcal K}).$$ As well
known any two isomorphisms from $\mathcal K$ to itself are
homotopic, hence $h_2\circ h_1^{-1}\sim id_{\mathcal K}$ and we
obtain (\ref{n_0}).
  Now to prove the commutativity of the left triangle of the diagram we have to prove that
\begin{equation}\label{n1}(g^A\otimes id_{\mathcal K})\circ (\rho^{A}\otimes id_{\mathcal K})\circ
(q\beta^A\otimes id_{\mathcal K}) \sim
 \gamma^{S^2A\otimes \mathcal K}.
  \end{equation} Like in Lemma \ref{K-theory} we will reduce the general case
  to the case $A = \mathbb C$ using the map $s^A: qA\to q\C\otimes
  A$ that was introduced in the proof of Lemma \ref{K-theory}.

  \noindent The right-hand side of (\ref{n1}) is  \begin{equation}\label{n2}\gamma^{S^2A\otimes \mathcal K} =
  (\gamma^{\C}\otimes id_{\sdAK})\circ (s^{\sdAK}\otimes id_{\mathcal
  K}), \end{equation} that is the diagram
  $$\xymatrix{q(S^2A\otimes K)\otimes
  K\ar[rr]^-{\gamma^{S^2A\otimes K}}\ar[dr]_-{s^{S^2A\otimes
  K}\otimes id_K} & & S^2A\otimes K\otimes K\\
  & q\C\otimes S^2A\otimes K\otimes K\ar[ur]_-{\gamma^{\C}\otimes
  id_{S^2A\otimes K}} &}$$
  commutes. It can be easily checked by comparing of the left-hand side and
the right-hand side of (\ref{n2}) on elements of $q(S^2A\otimes
K)\otimes K$ of the form $$\left(i_1^{\sdAK}(\phi\otimes a\otimes
S) - i_2^{\sdAK}(\phi\otimes
  a\otimes S)\right)\otimes T$$ and $$i_1^{\sdAK}(\phi_0\otimes a_0\otimes S_0)\left(i_1^{\sdAK}(\phi\otimes a\otimes S) -
  i_2^{\sdAK}(\phi\otimes
  a\otimes S)\right)\otimes T,$$ $\phi, \phi_0\in S^2\C$, $a, a_0\in A$, $T, S, S_0\in \mathcal
  K$, that span a dense subspace in $q(S^2A\otimes
K)\otimes K$.

  \noindent Clearly the left-hand side of (\ref{n1}) is equal to
  $\left(g^A\circ \rho^A\circ q\beta^A\right)\otimes id_{\mathcal K}$. We
  assert that \begin{equation}\label{n3} g^A\circ \rho^A\circ q\beta^A =
  (g^{\C}\otimes id_A)\circ (id_{q\C}\otimes \beta^{\C}\otimes
  id_A)\circ s^{\sdAK}, \end{equation} that is that the diagram
  $$\xymatrix{q(S^2A\otimes
  K)\ar[r]^-{q\beta^A}\ar[d]^{s^{S^2A\otimes K}}&q(A\otimes
  K)\ar[r]^-{\rho^A}&qA\otimes K\ar[d]^-{g^A}\\
  q\C\otimes S^2A\otimes K\ar[dr]_-{id_{q\C}\otimes
  \beta^{\C}\otimes id_A} & & S^2A\otimes K\\
  & q\C\otimes A\otimes K\ar[ur]_-{g^{\C}\otimes id_A}&}$$
  commutes. Indeed it is straightforward to show that the
  left-hand side and the right-hand side of (\ref{n3})
  asymptotically agree on elements of $q(S^2A\otimes K)$ of the
  form $$i_1^{\sdAK}(\phi\otimes a\otimes T) - i_2^{\sdAK}(\phi\otimes a\otimes
  T$$ and $$i_1^{\sdAK}(\phi_0\otimes a_0\otimes T_0)
  (i_1^{\sdAK}(\phi\otimes a\otimes T) - i_2^{\sdAK}(\phi\otimes a\otimes
  T)),$$  $\phi, \phi_0\in S^2\C$, $a, a_0\in A$, $T, T_0\in \mathcal
  K$, that span a dense subspace in $q(S^2A\otimes K)$.
Now,  by (\ref{n2}), (\ref{n3}), to get (\ref{n1})
  it remains to prove that $$(\gamma^{\C}\otimes id_{S^2\!A \!\otimes \!K})\circ
  (s^{S^2\!A \!\otimes \!K}\otimes id_{\mathcal K})\! \sim \!(g^{\C}\otimes id_{A\otimes
  K})\circ (id_{q\C}\otimes \beta^{\C}\otimes id_{A\otimes
  K})\circ (s^{S^2\!A \!\otimes \!K}\otimes id_{\mathcal K})$$ or, equivalently,
  $$\gamma^{\C}\otimes id_{S^2\C} \sim
  g^{\C}\otimes (id_{q\C}\otimes \beta^{\C}).$$ For that note that $\gamma^{\C}\otimes id_{S^2\C}$ and
  $g^{\C}\otimes (id_{q\C}\otimes \beta^{\C})$ induce the same
  homomorphisms in K-theory.
  Indeed they both send the generator $$\left(\left[\pr\right] -
  \left[p_{Bott}\right]\right)\otimes \left(\left[i_1^{\C}(1)\right] -
  \left[i_2^{\C}(1)\right]\right)$$ of $K_0(S^2\C\otimes q\C)$ to
  the generator $$\left[\pr\right] -
  \left[p_{Bott}\right]$$ of $K_0(S^2\C)$.

  \noindent This implies that $\gamma^{\C}\otimes id_{S^2\C}$ and
  $g^{\C}\otimes (id_{q\C}\otimes \beta^{\C})$ are homotopic
  because, as is well known, $$[[S^2\C\otimes
  q\C\otimes \mathcal K, S^2\C\otimes \mathcal K]]\cong Z \cong  Hom\left(K_0(S^2\C\otimes q\C),
  K_0(S^2\C)\right).$$ \qed

 Let $\psi\in G(B, A)$ and as before
$\phi:qB\to A\otimes \mathcal K$ be  an asymptotic morphism such
that $(id_A\otimes i)\circ(\phi\otimes id_{\mathcal K})\sim \psi$.
Define an asymptotic morphism $b(\psi)\in G(B, S^2A)$ by the
composition $$\xymatrix { qB\otimes K\ar[r]^{\theta_B}&
q^2B\otimes K\ar[rr]^-{q\phi\otimes id_K}&&q(A\otimes K)\otimes
K\ar[rr]^{\rho^A\otimes id_K}&&qA\otimes K\otimes K\\
&&& \ar[r]^-{id_{qA}\otimes i}&qA\otimes K\ar[r]^{g^A}&
S^2A\otimes K.}$$ Thus the map  $b: G(B, A)\to G(B, S^2A)$ is
defined by formula
$$ b(\psi) = g^A\circ(id_{qA}\otimes i)\circ(\rho^A\otimes id_{\mathcal K})\circ(q\phi\otimes id_{\mathcal K})\circ
\theta_B$$for any $\psi\in G(B, A)$.

\noindent Let $\beta^A_G: G(B, S^2A)\to G(B, A)$ be the map induced
by $\beta^A$.

\begin{prop}\label{GBott}  $b: G(B, A) \to G(B, S^2A)$ is a semigroup
isomorphism with inverse $\beta^A_G$.
\end{prop}
{\bf Proof.}  Obviously $b$ and $\beta^A_G$ are semigroup
homomorphisms so we have to check only the following:

\medskip

(i) $(\beta^A_G\circ b)(\psi) \sim \psi$ for any $\psi\in G(B,
A)$,

\medskip

(ii) $(b\circ \beta^A_G)(\psi) \sim \psi$ for any $\psi\in G(B,
S^2A)$.
\medskip

\noindent (i): Let $\psi\in G(B, A)$ and $\phi:qB\to A\otimes
\mathcal K$ be such an asymptotic morphism that $ \psi\sim
(id_A\otimes i)\circ(\phi\otimes id_{\mathcal K})$. Then
\begin{multline*} (\beta^A_G\circ b)(\psi) = \beta^A\circ g^A\circ
(id_{qA}\otimes i)\circ (\rho^A\otimes id_{\mathcal K})\circ
(q\phi\otimes
id_{\mathcal K})\circ \theta_B \stackrel{Lemma \ref{K-theory}}\sim \\
\gamma^A\circ (id_{qA}\otimes i)\circ (\rho^A\otimes id_{\mathcal
K})\circ
(q\phi\otimes id_{\mathcal K})\circ \theta_B \stackrel{Lemma \ref{1}} = \\
(id_A\otimes i)\circ \gamma^{A\otimes \mathcal K}\circ (q\phi\otimes
id_{\mathcal K})\circ \theta_B \stackrel{Lemma \ref{2}}=
(id_A\otimes i)\circ (\phi\otimes id_{\mathcal K})\circ
\gamma^{qB}\circ \theta_B \sim \\\psi\circ \gamma^{qB}\circ \theta_B
\stackrel{[1]} \sim \psi.
\end{multline*}

\noindent (ii): Let $\psi\in G(B, S^2A)$ and $\phi: qB\to
S^2A\otimes \mathcal K$ be such an asymptotic morphism that
$\psi\sim (id_{S^2A}\otimes i)\circ (\phi\otimes id_{\mathcal
K})$. Then
\begin{multline*} (b\circ \beta^A_G)(\psi) = g^A\circ
(id_{qA}\otimes i)\circ (\rho^A\otimes id_{\mathcal K})\circ
(q(\beta^A\circ \phi)\otimes id_{\mathcal K})\circ \theta_B = \\
g^A\circ (id_{qA}\otimes i)\circ (\rho^A\otimes id_{\mathcal
K})\circ (q\beta^A\otimes id_{\mathcal K})\circ (q\phi\otimes
id_{\mathcal K})\circ \theta_B \stackrel{Lemma \ref{new}}\sim
\\
(id_{S^2A}\otimes i)\circ \gamma^{S^2A\otimes \mathcal K}\circ
(q\phi\otimes id_{\mathcal K})\circ \theta_B \stackrel{Lemma
\ref{2}} = (id_{S^2A}\otimes i)\circ (\phi\otimes id_{\mathcal
K})\circ \gamma^{qB}\circ \theta_B \sim \\ \psi\circ
\gamma^{qB}\circ \theta_B \stackrel{[1]} {\sim} \psi.
\end{multline*}  \qed

\begin{prop}\label{G} $f^A_G\circ g^A_G = id$,

$\;\;\;\;\;\;\;\;\;\;\;\;\;\;\;\;\;\;\;\;\;\;\; g^A_G\circ f^A_G =
id$.
\end{prop}
Here $id$ means both the identity map from $G(B, S^2A)$ into
itself and the identity map from $G(B, qA)$ into itself.

 {\bf
Proof.} Consider the following diagram

$$
 \xymatrix{ & G(B, A) & \\
 G(B, S^2A)\ar[ur]^{\beta^A_G} \ar@<-1ex>[rr]_{f_E^A}& &G(B,
 qA) \ar[ul]_{\gamma_G^A} \ar@<-1ex>[ll]_{g_G^A}\\
 }
$$
 We shall prove that it commutes and this will imply the statement of the proposition.

\noindent By Propositions \ref{GBott} and \ref{Gq} $\beta^A_G$ and
$\gamma^A_G$ are isomorphisms. By Lemma \ref{K-theory}
$\beta^A_G\circ g^A_G = \gamma^A_G$ whence
\begin{equation}\label{old3}g^A_G = (\beta^A_G)^{-1}\circ \gamma^A_G
\end{equation} By Corollary \ref{sure} $\gamma^A_G\circ f^A_G =
\beta^A_G$ whence
\begin{equation}\label{old4}f^A_G = (\gamma^A_G)^{-1}\circ \beta^A_G
\end{equation} From (\ref{old3}) and (\ref{old4}) we obtain the assertions of the proposition. \qed
\subsection{Main result}
\begin{theorem}\label{main}
\item (i) $g^A\circ f^A \sim id_{\sdAK}$,

 \item (ii) $f^A\circ g^A \sim id_{\qAK}$.
\end{theorem}

{\bf Proof.} (i) By Proposition \ref{E} $g_E^A\circ f_E^A = id$
whence $g^A\circ f^A\circ \phi \sim \phi$ for any $\phi \in E(B,
\sdAK)$. Set $B = A\otimes \mathcal K$, $\phi = id_{\sdAK}$. Then
$$id_{\sdAK} \sim g^A\circ f^A\circ id_{\sdAK} = g^A\circ f^A.$$
  (ii) By Proposition \ref{G} $f^A_G\circ g^A_G = id$
whence $f^A\circ g^A \circ \phi \sim \phi$ for any $\phi \in
[[qB\otimes \mathcal K, \qAK]]$. Setting $B = A$, $\phi =
id_{\qAK}$ we get $$id_{\qAK} \sim f^A \circ g^A.$$
 \qed

\medskip

So $C^*$-algebras $\sdAK$ and $\qAK$ are asymptotically equivalent
 and we obtain immediately

\begin{cor}\label{E-theory}  $E(A, B) = [[qA, B\otimes \mathcal K]]$ for every $C^*$-algebras $A$ and $B$.
\end{cor}

\begin{cor}\label{Loring} Let $A$ be a nuclear $C^*$-algebra and $B$ be any $C^*$-algebra.
Then every asymptotic morphism from $qA$
 to $B\otimes \mathcal K$ is homotopic to a $*$-homomorphism from $qA$ to $B\otimes \mathcal K$.  \end{cor}

{\bf Proof.} Let $\phi_t\in [[qA, B\otimes \mathcal K]]$. Since
$A$ is nuclear  $I_{A, B}$ is an isomorphism (\cite{CH}). Define a
$\ast$-homomorphism $\psi_0:  \qAK \to B\otimes \mathcal K$ by
\begin{equation}\label{1end}\psi_0 = I_{A, B}^{-1}\left((id_B\otimes i)\circ
(\phi_t\otimes id_{\mathcal K})\circ f^A\right).\end{equation} By
\cite{Th} there exists a $\ast$-homomorphism $\psi: qA\to B\otimes
\mathcal K$ such that
\begin{equation}\label{2end}(id_B\otimes i)\circ (\psi\otimes id_{\mathcal K})\sim
\psi_0.\end{equation} We will prove that $\phi_t\sim \psi$. By
Theorem \ref{ind}
\begin{equation}\label{nonimp}I_{A, B}(\psi_0) = \psi_0\circ f^A.
\end{equation} By (\ref{1end}) the left-hand side of (\ref{nonimp}) is $I_{A,
B}(\psi_0) = (id_B\otimes i)\circ (\phi_t\otimes id_{\mathcal
K})\circ f^A$. By (\ref{2end}) the right-hand side of (\ref{nonimp})
is $\psi_0\circ f^A \sim (id_B\otimes i)\circ (\psi\otimes
id_{\mathcal K})\circ f^A$. So
$$(id_B\otimes i)\circ (\phi_t\otimes id_{\mathcal K})\circ f^A\sim (id_B\otimes i)\circ (\psi\otimes id_{\mathcal K})\circ
f^A,$$ $$(id_B\otimes i)\circ (\phi_t\otimes id_{\mathcal K})\circ
f^A\circ g^A\sim (id_B\otimes i)\circ (\psi\otimes id_{\mathcal
K})\circ f^A\circ g^A,$$ and by Theorem \ref{main} we obtain
$$(id_B\otimes i)\circ (\phi_t\otimes id_{\mathcal K})\sim (id_B\otimes i)\circ
(\psi\otimes id_{\mathcal K})$$ whence $\phi_t\sim \psi$.
 \qed

\medskip

{ \it Acknowledgements}. The author thanks V. M. Manuilov for
suggesting the problems treated in the paper and for many useful
discussions.

\end{document}